\documentclass[a4paper,11pt]{article}

\usepackage[latin1]{inputenc}
\usepackage{amsmath,amsfonts,amssymb}

\author{Jonas Lindstrøm Jensen}
\date{Department of Mathematics, University of Aarhus, Denmark}
\title{Correction of paper published in J. Combinatorial Theory 21, 1976: On the Existence of Hadamard Matrices}

\begin{document}
\maketitle

\section{Introduction}
An Hadamard Matrix is a square matrix with entries $\pm 1$, eg. a $n \times n$-matrix $H$ is an Hadamard Matrix if
\[ H H^\top = nI. \]
It is shown that for $n > 2$, it is a necessary condition that $n \equiv 4 \pmod n$, and it is conjectured that it is also a sufficient condition.
For a natural number $q>3$ it is proved in \cite{lit:seberry} that there exists Hadamard Matrices of order $2^s q$ for $s > t$, where $t \leq [2\log_2(q-3)]$. In this paper I will show that the bound on $t$ is wrong. It is true that there is some $t$ to make the claim true, but the bound on $t$ is wrong.

\section{Structure of the proof}
In \cite{lit:seberry} a corollary to a two dimensional version of Frobenius Coin Problem is used. The theorem states that for relatively prime integers $x,y$, any integer $N > (x-1)(y-1)$ can be written in the form $ax + by$ for some nonnegative integers $a,b$. The corollary states that given $x = v+1$ and $y=v-3$ where $v \geq 9$ is odd, there exist nonnegative integers $a,b$ such that
\[ a(v+1) + b(v-3) = 2^t, \]
for some $t$. The proof of Corollary 7 goes like this (I have filled in some details that is left out in \cite{lit:seberry}).

Let $g = \gcd(v+1, v-3)$, then $g \in \{1,2,4\}$, and hence $g=2^d$ for some $d$. Let
\[ N = \left( \frac{v+1}{g} - 1 \right) \left( \frac{v-3}{g} - 1 \right), \]
and let $2^k$ be the smallest power of 2 greater than N. By the theorem we have nonnegative $a,b$ such that
\[ 2^k = a\frac{v+1}{g} + \frac{v-3}{g}, \]
and since $g = 2^d$ we have
\[ 2^t = a(v+1) + b(v-3), \]
where $t=k+d$.

In Lemma 9 this result is used to show that there exist nonnegative $a,b$ such that $a(v+1) + b(v-3) = 2^t$ to construct an Hadamard Matrix of order $2^{t+1} v$. This only holds for $v$ being a prime $\equiv 1 \pmod 4$, and in that case $g=2$ in the proof of the corollary and  hence $t = k+1$. A similar proof is done for $v \equiv 3 \pmod 4$ where $g=4$ and $t=k+2$. Since the Kronecker product of two Hadamard Matrices, the prime factorization of $q$ can be used to construct an Hadamard Matrix of order $2^s q$ where $s \geq t$ for a sufficiently large $t$.

\section{The error}
The error comes when trying to estimate how big $t$ has to be. In \cite{lit:seberry} it is stated that $t < [2 \log_2(q - 3)]$ is enough for each prime factor, and since
\begin{equation}\label{eq:multest}
 [2 \log_2(p-3)] + [2 \log_2(q-3)] < [2\log2(pq-3)],
\end{equation}
that bound is preserved under multiplication, and hence when used on all prime factors.

But this is not the case, since for $v \equiv 1 \pmod 4$ has $g=2$ and choose $k$ such that
\begin{equation}\label{eq:kest}
 2^k > \left( \frac{v+1}{2} - 1 \right) \left( \frac{v-3}{2} - 1 \right) = \frac{v-1}{2} \frac{v-5}{2}, \
\end{equation}
which implies
\[ k > \log_2\left( \frac{1}{4}(v-1)(v-5) \right) = \log_2(v-1) + \log_2(v-5) - 2. \]
So if we choose $k = [2\log_2(v-3)] - 1$, we can ensure (\ref{eq:kest}) to be true. But as we stated earlier, this gives us relation $a(v+1) + b(v-3) = 2^t$ for $t=k+1$, and Lemma 9 uses this relation to give us an Hadamard Matrix of order $2^{t+1} v$, which is double the bound stated in \cite{lit:seberry}. With the new bound $t < [2 \log_2(q - 3)] + 1$, it is not possible to make an estimate like (\ref{eq:multest}), and the estimate cannot be extended from the prime factors to the product $q$.

\section{A numerical example}
To make things a bit clearer, I provide a numerical example. Let $v=17$. It is claimed that Lemma 9 gives an Hadamard Matrix of order $2^{t+1} v$ where $t = [2 \log_2(v-3)]-1$. So in our case $t=6$. In the proof of Corollary 7, we choose $k$ such that $2^k$ is greater than
\[ N = \left( \frac{17+1}{2} - 1 \right) \left( \frac{17-3}{2} - 1 \right) = 8 \cdot 6 = 48, \]
hence $k=6$. The theorem gives us a relation
\[ a \frac{17+1}{2} + b \frac{17-3}{2} = 2^6, \]
and by multiplication with 2 we get
\[ 18a + 14b = 2^7. \]
Lemma 9 in \cite{lit:seberry} uses this to provide us with an Hadamard Matrix of order $2^8 \cdot 13$, but the bound on the exponent given in \cite{lit:seberry} is $t \leq [2 \log_2(14)] = 7$.


\begin{thebibliography}{WWW} 
\bibitem{lit:seberry} J. Seberry Wallis (1976), On the Existence of Hadamard Matrices, \emph{J. Comb. Theory (A) \textbf{21}}, 188-195
\end{thebibliography}
\end{document}